\documentclass{amsart}
\usepackage[all]{xy}
\usepackage{amsmath}
\usepackage{latexsym}
\RequirePackage{amssymb}

\newtheorem{prop}[subsection]{Proposition}
\newtheorem{cor}[subsection]{Corollary}
\newtheorem{thm}[subsection]{Theorem}
\newtheorem{lem}[subsection]{Lemma}
\newcommand{\g}{\mathfrak{g}}
\newcommand{\h}{\mathfrak{h}}

\newcommand{\plie}{\bf{p\textrm{-} Lie}}
\newcommand{\pleib}{\bf{p\textrm{-} Leib}}

\newcommand{\lp}{\dashv}
\newcommand{\rp}{\vdash}

\newcommand{\Proof}{\emph{Proof.\;\;}}
\begin{document}

\title{On restricted Leibniz algebras}

\bigskip
\bigskip

% Information for first author
\author{Ioannis Dokas}
\address{Institut de Recherche Mathematique Avanc\'{e}e, Universit\'{e}
Louis Pasteur, 7 rue Ren\'{e}-Descartes, 67084 Strasbourg, France}
\email{dokas@math.u-strasbg.fr}

%Information for second author
\author{Jean-Louis Loday}
%Address of record for the research reported here
%\address{}
%Current address
%\curraddr{}
\address{Institut de Recherche Mathematique Avanc\'{e}e, Universit\'{e}
Louis Pasteur, 7 rue Ren\'{e}-Descartes, 67084 Strasbourg, France}
\email{loday@math.u-strasbg.fr}
%\thanks will become a 1st page footnote.
%\thanks{}

%    General info

%\date{ }

%\dedicatory{}


\begin{abstract}
In this paper we prove that in prime characteristic there is a
functor $-_{\pleib}$ from the category of diassociative algebras
to the category of restricted Leibniz algebras, generalizing the functor from associative algebras to restricted Lie algebras. Moreover we
define the notion of restricted universal enveloping diassociative
algebra $Ud_{p}(\g)$ of a restricted Leibniz algebra $\g$ and we
show that $Ud_{p}$ is left adjoint to the functor
$-_{\pleib}$. We also construct the 
 restricted enveloping  algebra, which classifies the restricted Leibniz modules.
 In the last section we put a restricted pre-Lie structure on the tensor product of a Leibniz algebra by a Zinbiel algebra.
\end{abstract}

\maketitle

%\tableofcontents

\section*{Introduction}

In characteristic $p$ the notion of Lie algebra is fruitfully replaced by the notion of \emph{restricted Lie algebra}, that is a Lie algebra equipped with a formal $p$-th power called the Frobenius map. A non-antisymmetric version of Lie algebras has been introduced by the second author in \cite{L1} under the name \emph{Leibniz algebras}. In \cite{DZ}  Dzhumadil'daev and S.A. Abdykassymova ~introduced the notion of \emph{restricted Leibniz algebra} in characteristic $p$.

It is well-known that an associative algebra gives rise to a restricted Lie algebra by taking the $p$-th power as Frobenius map. In the Leibniz case, we prove that the Leibniz algebra associated to a diassociative algebra is in fact a restricted Leibniz algebra.

 In the Lie case, the functor from associative algebras to restricted Lie algebras admits a left adjoint $U_p$ and it is known that,  for a restricted Lie algebra $\g$, the associative algebra $U_p(\g)$ classifies the restricted Lie modules over $\g$. In the Leibniz case, these two roles of $U_p$ are played by two different objects. First, we prove that  the restricted enveloping diassociative algebra $Ud_{p}(\g)$ provides a functor which is left adjoint to the forgetful functor. Second, we show that the 
 restricted Leibniz $\g$-modules are classified by an associative algebra $UL_{p}(\g)$
 that we construct explicitly.

 In the last part we show that the tensor product of a Leibniz algebra and a Zinbiel algebra (dual notion for Koszul duality), which is known to be a Lie algebra, has a finer structure: it is a \emph{pre-Lie algebra}. This result is valid in any characteristic. In characteristic $p$ it turns out that this pre-Lie algebra is in fact a restricted pre-Lie algebra in the sense of  Dzhumadil'daev, cf.~\cite{DZ2}.

We thank Bruno Vallette for a careful reading of a first version of this paper.
\bigskip

In this paper we denote by $k$ a field of prime characteristic
$p$ except at the beginning of section \ref{Lie} where no characteristic hypothesis is made.

\section{Diassociative algebras and restricted
Leibniz algebras.}

A Leibniz algebra is a non-commutative version of Lie algebra. Thus, before we proceed with Leibniz algebras let us recall some facts in the Lie framework which led N. Jacobson to introduce the notion of restricted Lie algebra. In many cases in which Lie algebras arise naturally in prime characteristic one finds structures richer than the ordinary Lie algebras. For example, any associative algebra $A$ over $k$ gives rise to a pair $(A_{Lie},[p])$ where $A_{Lie}$ is the Lie algebra with bracket given by $[a,b]:=ab-ba$ where $a,b\in A$ and $[p]: A_{Lie} \rightarrow A_{Lie}$ is the `Frobenius mapping' $a\mapsto a^{p}$. Another example is obtained by considering the Lie algebra $Der(\frak{U})$ of derivations of a not necessarily associative algebra $\frak{U}$ over $k$. It can be proved (see \cite{Jac}) that $Der(\frak{U})$ is closed under the map $[p]:D \mapsto D^{p}$ where $D\in Der(\frak{U})$ and the map $[p]$ verifies specific relations. These facts lead to the following definition.

\subsection{Definition \cite{Jac}}
A \emph{restricted Lie algebra} is a pair $(L,[p])$ where $L$ is a Lie algebra over $k$ and $[p]: L \rightarrow L, \; x\mapsto x^{[p]}$ is map which satisfies the following relations:
\begin{eqnarray*}
(\alpha x)^{[p]} &=& \alpha^{p}x^{[p]}, \;\; \alpha \in k,\; x\in L \\
\,[x,y^{[p]}\,]   &=& \,[\cdots \,[x,\underbrace{y\,],y\,]\cdots,y}_{p}\,], \;\;x,y \in L\\
(x+y)^{[p]}       &=& x^{[p]}+ y^{[p]}+ \sum_{i=1}^{i=p-1}s_{i}(x,y),
\end{eqnarray*}
where $is_{i}(x,y)$ is the coefficient of $\lambda^{i-1}$ in $\,[\cdots \,[x,\underbrace{(\lambda x+y)\,],(\lambda x+y)\,]\cdots,(\lambda x+y)}_{p-1}\,]$.

A \emph{restricted} morphism $f: L_{1} \rightarrow L_{2}$ is a Lie morphism such that $f(x^{[p]})=f(x)^{[p]}$. We denote by $\plie$ the category of restricted Lie algebras.

\subsection{Example \cite{Jac}}\label{examplerestricted}
Any associative algebra $A$ over $k$ has the structure of restricted algebra with $p$-map given by $a \mapsto a^{p}$. In particular, from the Leibniz rule we obtain that the derivation algebra $Der(\frak{U}) \subset gl(\frak{U})$ of a not necessarily associative algebra $\frak{U}$ is a restricted Lie algebra.

Let $L$ be a restricted Lie algebra, and let $U(L)$ be the universal enveloping algebra of $L$. The \emph{universal restricted enveloping algebra} is defined by $U_{p}(L):= U(L)/\{x^{p}-x^{[p]}\;\; |\;x\in L\}$. A module over $U_{p}(L)$ is called  \emph{restricted} Lie module.

From example \ref{examplerestricted} above we see that there is a functor $\bf{Ass} \rightarrow \bf{p \textrm{-}Lie}$. The left adjoint of this functor is the restricted enveloping functor $U_{p}(L)$. A natural question which arises is wether there is an analogous statement in the Leibniz framework in prime characteristic. Let us now recall the notion of Leibniz algebra, for more
details about Leibniz algebras the reader could consult \cite{L1},\cite{LP}.

\subsection{Definition}
A $\emph{Leibniz algebra}$ over $k$ is a $k$-module $\g$ equipped
with a bilinear map, called $\emph{bracket}$,
\[\,[-,-\,]: \g \times \g\rightarrow \g\] satisfying the
$\mathit{Leibniz\; identity:}$
\[\,[x,\,[y,z\,]\,]=\,[\,[x,y\,],z\,]-\,[\,[x,z\,],y\,],\] for all
$x,y,z\in \g$

We denote by $\bf{Leib}$ the category of Leibniz algebras over $k$.

\subsection{Example}
Any Lie algebra is a Leibniz algebra.

When we replace Lie algebras by Leibniz algebras the role of
associative algebras is played by the diassociative algebras
(see \cite{L3}). We recall the definition:

\subsection{Definitions}\label{diassociative}
A $\emph{diassociative}$ algebra is a $k$-module $D$ equipped with
two $k$-linear maps $\lp,\; \rp : D \otimes D \rightarrow D$
called respectively the $\emph{left\; product}$ and the
$\emph{right product}$ such that the products $\lp$ and $\rp$ are
associative and satisfy the following relations:
\begin{eqnarray}
x \lp (y \lp z) &=& x \lp (y\rp z), \\
(x \rp y) \lp z &=& x \rp (y \lp z), \\
(x \lp y) \rp z &=& (x \rp y )\rp z.
\end{eqnarray}

Let $D$ and $D'$ be diassociative algebras. A \emph{morphism} of
diassociative algebras $f: D \rightarrow D'$ is a $k$-linear map
such that:
\[f(x \lp y)=f(x) \lp f(y) \quad and \quad f(x \rp y)=f(x) \rp f(y)\]
for all $x,y \in D$.

The category of diassociative algebras over $k$ is denoted by $\bf{Dias}$.

\begin{prop}
Let $D$ be a diassociative algebra. Then $D$ endowed with the following bracket
\[ \,[x,y\,]:=x\lp y -y \rp x \] is a Leibniz algebra, denoted by
$D_{Leib}$.
\end{prop}

\Proof See proposition 4.2 in \cite{L3}.

\medskip

Here we would like to note that in the Lie context the solution of many problems in prime characteristic requires only the knowledge of whether a $p$-map can be introduced in a given Lie algebra $L$ over $k$. This led to the concept of \emph{restrictable} Lie algebras (see \cite{Far}). In particular a Lie algebra $L$ over $k$ is called \emph{restrictable} if $(\textrm{Ad }x)^{p}$ is an  inner derivation for all $x\in L$ (here $\textrm{Ad }x (y) = [y,x]$). In other words for all $x\in L$ there exist $x^{[p]} \in L$ such that $(\textrm{Ad }x)^{p}=(\textrm{Ad }x^{[p]})$. From Theorem $11$ in \cite{Jac} we obtain that a Lie algebra $L$ is restrictable if and only if there is a $p$-map $[p]: L\rightarrow L$ which makes $L$ a restricted Lie algebra. 

 In \cite{DZ} A.S Dzhumadil'daev and S.A. Abdykassymova introduce an analogue of the notion of \emph{restrictable} Lie algebra in the Leibniz framework. They call it  restricted Leibniz algebra, and we keep their terminology.

\subsection{Notation}
Let $\g$ be a Leibniz algebra, we denote by $r_{x}: \g
\rightarrow \g$ the right multiplication operator given by
$r_{x}(y):=\,[y,x\,]$ for all $x,y \in \g$.

\subsection{Definition \cite{DZ}}
A \emph{restricted Leibniz algebra} is a pair $(\g,[p])$ where
$\g$ is a Leibniz algebra $\g$ over $k$, and $[p]: \g \rightarrow
\g,\; x \mapsto x^{[p]}$ is a map (called the $p$-map) such that
$r_{x}^{p}=r_{x^{[p]}}$ for all $x\in \g$.

Let $(\g,[p])$ and $(\h,[p])$ be restricted Leibniz
algebras. A \emph{restricted morphism} $f: \g \rightarrow \h$ is a
Leibniz morphism such that $f(x^{[p]})=f(x)^{[p]}$ for
all $x\in \g$. We denote the category of restricted Leibniz algebras over $k$ by $\pleib$.

\subsection{Example}
Any restricted Lie algebra is restricted as a Leibniz algebra.

The next theorem shows that 
the associated Leibniz algebra $D_{Leib}$ of a diassociative
algebra $D$ over $k$ has the structure of restricted Leibniz
algebra. We first prove a technical lemma.

\begin{lem}\label{lemdias}
If $D$ is  a diassociative algebra then we have \[x\lp
(\cdots(\underbrace{y\lp y)\lp y) \cdots)\lp y}_{n})=x\lp
(\cdots(\underbrace{y\rp y)\rp y) \cdots)\rp y}_{n})\] for all
$x,y \in D$ and $n\geq 1$. We adopt the notation $y^{\lp\;
n}:=(\cdots(\underbrace{y\lp y)\lp y) \cdots)\lp y}_{n}$ and 
$y^{\rp\;n }:=(\cdots(\underbrace{y\rp y)\rp y) \cdots)\rp
y}_{n}$.
\end{lem}

\Proof In order to prove the identity we proceed by induction on
$n$. For $n=2$ this identity follows from axiom  $(1)$ of the definition
\ref{diassociative}. We suppose that the relation is true for $n$. Then $x\lp
y^{\lp\; (n+1)}=x\lp (y^{\lp\; n}\lp y)=(x \lp y^{\lp\; n})\lp y$.
By induction hypothesis we get $(x \lp y^{\lp\; n})\lp y=(x\lp
y^{\rp\; n})\lp y=x\lp (y^{\rp\; n}\lp y)$. Finally, by axiom $(1)$
of the definition \ref{diassociative}, we obtain that $x\lp (y^{\rp\; n}\rp
y)=x\lp y^{\rp\; (n+1)}$. \hfill $\square$

\begin{thm}\label{diasso}
Let $D$ be a diassociative algebra over $k$. If $[p]: D
\rightarrow D$ is the map given by $x \rightarrow x^{\rp \; p}$,
then $(D_{Leib}, [p])$ is a restricted Leibniz algebra.
\end{thm}

\Proof For all $y \in D$ we will denote by $R_{y}^{\lp},
L_{y}^{\rp}: D\rightarrow D$ the maps defined respectively  by
$R_{y}^{\lp}(x):=x\lp y$ and $L_{y}^{\rp}(x):=y\rp x$. The Leibniz
bracket, with this notation, is given by
$\,[x,y\,]=(R_{y}^{\lp}-L_{y}^{\rp})(x)$\, for all\, $x,y \in D$.
Moreover, from relation $(2)$ of the definition \ref{diassociative} we have that
$R_{y}^{\lp}L_{y}^{\rp}=L_{y}^{\rp}R_{y}^{\lp}$. Therefore by the
binomial formula  we obtain that
$(R_{y}^{\lp}-L_{y}^{\rp})^{n}=\sum_{i=0}^{i=n}\binom{n}{i}(R_{y}^{\lp})^{i}(-L_{y}^{\rp})^{n-i}$.

Besides in prime characteristic $p$ we have $\binom{p}{i}=0 \ 
(mod\;p)$ \,if\, $0<i<p$. Thus, for $p=n$ we get that
$(R_{y}^{\lp}-L_{y}^{\rp})^{p}=(R_{y}^{\lp})^{p}-(L_{y}^{\rp})^{p}$.
Therefore we have:
\begin{eqnarray*}
\,[\cdots\,[x,\underbrace{y\,],y\,]\cdots y}_{p}\,]
&=&(R_{y}^{\lp})^{p}(x)-(L_{y}^{\rp})^{p}(x) \\
                                                    &=&(\cdots(x\lp
\underbrace{y)\lp y)\cdots )\lp
y}_{p})-(\underbrace{y\rp\cdots(y\rp(y}_{p}\rp x)\cdots ) \\
                                                    &=& x \lp (\cdots
(\underbrace{(y \lp y)\lp y)\cdots )\lp y}_{p})-(\underbrace{y\rp(
\cdots (y \rp (y \rp y}_{p})\cdots) \rp x \\
&=& x \lp y^{\lp \; p}-y^{\rp \;p} \rp x.
\end{eqnarray*}
Finally, by lemma \ref{lemdias} above we obtain:
\[\,[\cdots\,[x,\underbrace{y\,],y\,]\cdots y}_{p}\,]=x \lp y^{\rp \;
p}-y^{\rp\; p} \rp x=\,[x,y^{\rp\;p}\,].\]
\hfill $\square$

\subsection{Remark} Here we made of a choice to define a $p$-map. In \cite{Goi} Goichot introduced the quotient $D_{as}$ of $D$ by the relation $\lp\ =\ \rp$. In fact all the elements
$x\rp \cdots \rp x \lp \cdots \lp x$ have the same image in $D_{as}$ and each of them can be taken as $p$-map.

\begin{prop} Let $D$ be a diassociative algebra over $k$. The following formula holds in $D_{p\textrm{-}Leib}$ 
$$ [z, (x+y)^{[p]}]       = [z, x^{[p]}]+[z,  y^{[p]}]+ [z, \sum_{i=1}^{i=p-1}s_{i}(x,y)],$$
where the bracket involved in $s_{i}(x,y)$ is the Leibniz bracket  $x\lp y - y \rp x$.
\end{prop}

\Proof Since $\rp$ is an associative product we have the Jacobson formula 
$$(x+y)^{[p]}       = x^{[p]}+ y^{[p]}+ \sum_{i=1}^{i=p-1}s_{i}(x,y),$$
where the bracket used in $s_i$ is the Lie bracket $x\rp y - y \rp x$. Under left bracketing with $z$ we get the same element by replacing this bracket by the Leibniz bracket because of the relations $z\lp (x\lp y) = z\lp (x\rp y)$ and  $(x\lp y)\rp z =  (x\rp y)\rp z$. 
\hfill $\square$\\

From Theorem \ref{diasso} above we see that, in prime
characteristic, the structure of restricted Leibniz algebra arises
in a natural way from the structure of diassociative
algebra.  Moreover we have the following Proposition:

\begin{prop}
The following diagram of categories of algebras is commutative.
\begin{displaymath}
\xymatrix{ 
\bf{As} \ar[d] \ar[r] & 
         \plie \ar[d] \\
\bf{Dias} \ar[r] &  \pleib   }
\end{displaymath}
\end{prop}

\subsection{Example}
Let $D$ be an diassociative algebra over $k$, we denote by
$\mathcal{M}(D)$ the diassociative algebra of $n\times n$-matrices
with entries in $D$. Then
$\mathfrak{g}l_{n}(D):=\mathcal{M}(D)_{Leib}$ is a restricted
Leibniz algebra. Here we would like to mention that in zero
characteristic the homology of $\mathfrak{g}l_{n}(D)$ when $D$ has
a bar-unit, has been computed by Frabetti (see \cite{fr}).

\subsection{Example}
Let $A$ be an associative $k$-algebra equipped with a $k$-module
map $D: A\rightarrow A$ satisfying the condition:
\[D(a(Db))=DaDb=D((Da)b)\] for all\; $a,b \in A$. We define a bilinear
map on $A$ by \[\,[a,b\,]:=a(Db)-(Db)a.\] If we set $a \lp b:=
aDb$ and $b\rp a:= (Db)a$ for all $a,b \in A$, then we can easily
check that $(A,\lp,\rp)$ is a diassociative algebra. Therefore by
theorem \ref{diasso} above $A$ has the structure of restricted
Leibniz algebra with  map $[p]: A \rightarrow A$ given by $a
\mapsto a^{\rp \;p}=(Da)^{p-1}a$ for all $a \in A$.

In case that $D=id$, we notice that $(A, \,[-,-\,])$ is a
restricted Lie algebra with map given by $a \mapsto a^{p}$ for all
$a \in A$.

\section{Restricted universal enveloping diassociative algebra.}
\bigskip

Let us recall that the free diassociative algebra over the vector space $V$ is of the form $T(V)\otimes V \otimes T(V)$, cf. \cite{L3}.

The functor $(-)_{Leib}:\; \bf{Dias} \rightarrow \bf{Leib}$ has a
left adjoint functor $Ud : \bf{Leib} \rightarrow \bf{Dias}$ given
by
\[Ud(\g):=T(\g)\otimes \g \otimes T(\g)/\{[x,y]-x\lp y+y\rp x\;\;|x,y\in\g\}.\]
By theorem \ref{diasso} we obtain a functor $(-)_{p\textrm{-}Leib}:\; \bf{Dias}
\rightarrow \pleib$. In the following we construct a left
adjoint to the functor $(-)_{p\textrm{-}Leib}.$

\subsection{Definition}\label{uea}
Let $(\g,[p])$ be a restricted Leibniz algebra. Then the
\emph{restricted universal diassociative algebra} $Ud_{p}(\g)$ is
defined by:
\[T(\g)\otimes \g \otimes T(\g)/\big\{[x,y]-x\lp y+y\rp x,\; x^{[p]}-x^{\rp\;p}\;\;|\ x,y\in\g\big\},\]
where, in $T(\g)\otimes \g \otimes T(\g)$, we have made the following identification:
\begin{eqnarray*}
[x,y] & = & 1\otimes [x,y]\otimes 1\ ,\\
x\lp y& = &  1\otimes x\otimes y\ ,\\
y\rp x & = & y\otimes x \otimes 1\ , \\
 x^{[p]} & =  & 1 \otimes x^{[p]}\otimes 1 \ , \\
x^{\rp\;p} & = & \underbrace{x\otimes x \cdots \otimes x}_{p} \otimes 1\ .\\
\end{eqnarray*}

\begin{prop}
 Let $(\g,[p])$ be a restricted Leibniz algebra and $D$ a diassociative
algebra over $k$. Then we have the following natural bijection of sets.
\[Hom_{\bf{Dias}}(Ud_{p}(\g),D)\simeq Hom_{\pleib}(\g,D_{p\textrm{-}Leib}).\]
In other words the functor $Ud_{p}$ is left adjoint to the functor $(-)_{p\textrm{-}Leib}$.
\end{prop}

\Proof Let $f: \g \rightarrow D_{\pleib}$ be a morphism of
restricted Leibniz algebras. Since $T(\g)\otimes \g \otimes T(\g)$
is the free diassociative algebra on $\g$ there is a unique
extension of $f$ to a morphism $\hat{f}:T(\g)\otimes \g \otimes
T(\g) \rightarrow D$  of diassociative algebras. Moreover, we have
$\hat{f}([x,y]-x\lp y+y\rp x)=0$ and $\hat{f}(x^{[p]}-x^{\rp\;p})$
because $f$ is a restricted Leibniz morphism. Therefore, $\hat{f}$
induces a diassociative algebra morphism $\phi_{f}:Ud_{p}(\g)
\rightarrow D$.

Conversely, if $\phi: Ud_{p}(\g) \rightarrow D$ is a diassociative
algebra morphism then the restriction $f_{\phi}: k\otimes \g
\otimes k \rightarrow D$ of $\phi$ to $k\otimes \g \otimes k$ is a
restricted Leibniz morphism.

It is easy to check that $f_{(\phi_{f})}=f$ and
$\phi_{(f_{\phi})}=\phi$ and these two constructions give rise to
a bijection
\[Hom_{\bf{Dias}}(Ud_{p}(\g),D)\simeq Hom_{\bf{p-Leib}}(\g,D_{p\textrm{-}Leib}).\]
\hfill $\square$

\section{Restricted modules and  extensions.}

In order to compare with the Lie context let us recall that a strongly abelian  extension of restricted Lie algebra $\g$ by a restricted Lie module $M$ is an exact sequence of restricted Lie algebras:
\[ 0\rightarrow M \rightarrow \mathfrak{b}
\rightarrow \g \rightarrow 0 \] such that $[M,M]=0$ and
$M^{[p]}=0$.
\subsection{Remark}
There is a bijection between the second Hochschild cohomology group $H^{2}(\g ,M)$ (see \cite{Ho}) and the set of equivalence classes of strong abelian extensions. We note also that we can characterize extensions for which the $p$-map is not $0$ using the general scheme of Quillen-Barr-Beck's cohomology theory for details the reader may consult \cite{dok}.

By analogy we give the following definition:

\subsection{Definition}
An  \emph{abelian extension of restricted Leibniz algebras}
(extension of $\g$ by $M$) is an exact sequence of restricted
Leibniz algebras: \[(\mathcal{E}):\;\;\; 0\rightarrow M \rightarrow \mathfrak{b}
\rightarrow \g \rightarrow 0 \] such that $[M,M]=0$ and
$M^{[p]}=0$

Here we recall the definition of \emph{restricted} Leibniz module given by A. S. Dzhumadil'daev and S. A. Abdykassymova in \cite{DZ}.

\subsection{Definition \cite{DZ}}
 A \emph{restricted} Leibniz $\g$-module is a Leibniz $\g$-module such
that \[\,[m,x^{[p]}\,]=\,[\cdots [m, \underbrace{x],x]\cdots
,x}_{p}].\]

\begin{lem}
If $(\mathcal{E})$ is an extension of the restricted Leibniz algebra $\g$ by $M$
\[(\mathcal{E}):\;\; 0\rightarrow M \rightarrow \mathfrak{b}
\rightarrow \g \rightarrow 0\] with $[M,M]=0$ and $M^{[p]}=0$, then the vector space $M$ inherits the structure of restricted Leibniz module.
\end{lem}
\Proof Given such an extension we see that $M$ is a
$k$-module equipped with two actions (left and right) of $\g$
\[
 [-,-]: \g \times M \rightarrow M \quad  and \quad [-,-]:
M \times \g \rightarrow M
\] such that the following relations hold:
\begin{eqnarray}
\,[m,\,[x,y\,]\,] &=& \,[\,[m,x\,],y\,]-\,[\,[m,y\,],x\,], \\
\,[x,\,[m,y\,]\,] &=& \,[\,[x,m\,],y\,]-\,[\,[x,y\,],m\,],\\
\,[x,\,[y,m\,]\,] &=& \,[\,[x,y\,],m\,]-\,[\,[x,m\,],y\,],   \\
\,[m,x^{[p]}\,]   &=& \,[\cdots [m, \underbrace{x],x]\cdots
,x}_{p} \,],
 \end{eqnarray}
for all $x,y\in \g$ and $m\in M$.
Therefore $M$ is a restricted Leibniz module. \hfill $\square$

\subsection{Restricted universal enveloping algebra of a restricted
Leibniz algebra}
\bigskip

Let $\g^{l}$ and $\g^{r}$ be two copies of the Leibniz algebra
$\g$. We denote by $l_{x}$ and $r_{x}$ the elements of $\g^{l}$
and $\g^{r}$ corresponding to $x \in \g$. Let $T(\g^{l}\oplus
\g^{r})$ be the tensor k-algebra, which is associative and unital.

\subsection{Definition}\label{ruea}
Let $\g$ be a restricted Leibniz algebra. The \emph{restricted
universal enveloping algebra} $UL_{p}(\g)$ of $\g$ is defined by
$UL_{p}(\g):=T(\g^{l}\oplus \g^{r})/I_{p}$ where $I_{p}$ is the
two-sided ideal corresponding to the relations:
\begin{eqnarray}
r_{[x,y]}- r_{x}r_{y}-r_{y}r_{x}& =& 0,\\
l_{[x,y]}-l_{x}r_{y}-r_{y} l_{x}& =& 0, \\
(r_{y}+l_{y})l_{x} &= & 0, \\
r_{x^{[p]}}-r_{x}^{p} &= & 0,
\end{eqnarray}
for all $x,y \in \g$.

\begin{thm}
Let $\g$ be a restricted Leibniz algebra. The category of
restricted Leibniz $\g$-modules is equivalent to the category of
right modules over $UL_{p}(\g)$.
\end{thm}

\Proof Let $M$ be a restricted $\g$-module. Then we define actions
of $\g^{l}$ and $\g^{r}$ respectively by:
\[ m \cdot l_{x}:=[x,m], \quad and \quad m \cdot r_{x}:=[m,x].\]
for all $m \in M$ and $x,y \in \g$. These actions are extended by
composition and linearity to an action of $T(\g^{l}\oplus
\g^{r})$. Since $M$ is a restricted $\g$-module it is in particular a
Leibniz $\g$-module, therefore the elements of type $(8),(9),(10)$
act trivially. Besides, from axiom $(7)$ elements of type $(11)$
act trivially as well. Therefore $M$ is endowed with the structure
of right $UL_{p}(\g)$-module.

Conversely, let $M$ be a right $UL_{p}(\g)$-module. Then the
restriction of the action on $\g^{l}$ and $\g^{r}$ makes $M$ into a
Leibniz $\g$-module. Moreover, from relation $(11)$ of the
definition \ref{ruea} we obtain that $M$ is actually a restricted
Leibniz $\g$-module. \hfill $\square$

\section{Lie property of Koszul duality for Leibniz algebras}\label{Lie}

Let us recall two notions of algebras which are closely related to Lie 
and Leibniz algebras.

By definition a \emph{pre-Lie algebra} (also called Vinberg algebra, and right-symmetric algebra) is a vector space $A$ equipped 
with a binary operation denoted $\{-,-\}$ whose associator is 
right-symmetric, that is, satisfies the relation
$$\{\{x,y\},z\}- \{x,\{y,z\}\}= \{\{x,z\},y\}- \{x,\{z,y\}\}.$$
It is well-known that the antisymmetrized operation $[x,y]:= 
\{x,y\}-\{y,x\}$ is a Lie bracket. In other words, pre-Lie algebras are Lie admissible algebras.
\subsection{Example} Any associative algebra is a pre-Lie algebra.

\subsection{Notation}
Let $A$ be a pre-Lie algebra. We denote by $R_{a}: A \rightarrow A$ the right-multiplication operator given by $R_{a}(b):=\{b,a\}$ for all $a,b \in A$. Moreover, we denote the $n$-th power, by $a^{\{n\}}:=\{\cdots \{\{\underbrace{a,a\},a\}\cdots\},a}_{n}\}$ where $a\in A$ and $n\geq 1$.

Next, let us recall the notion of Zinbiel algebra which is the dual of Leibniz algebra in Koszul sense.

By definition a \emph{Zinbiel algebra} is a vector space $\mathfrak{R}$ equipped 
with a binary operation denoted $\prec$ which satisfies the relation
$$(a\prec b)\prec c = a\prec (b\prec c)+ a\prec (c\prec b) .$$

The following lemma is useful for computations.
\begin{lem}\label{zinb}
Let $\mathfrak{R}$ be a Zinbiel algebra over $k$. In any characteristic we have, for any $a,b\in \mathfrak{R}$:
$$ (\cdots((a\prec \underbrace{b)\prec b) \cdots )\prec b}_{n} = n!\  a\prec \underbrace{(b\prec (b\prec(  \cdots \prec b)))}_{n}.$$
In particular, in characteristic $p$, we get:
\[ (\cdots((a\prec \underbrace{b)\prec b) \cdots )\prec b}_{p}=0\] for all $a,b \in \mathfrak{R}$.
\end{lem}
\proof 
Let us denote $b^{\prec n}:= b\prec b^{\prec n-1}$ where $b^{\prec 1}=b$. We prove by induction that $b^{\prec n}\prec b= n  \ b^{\prec n+1}$. It is true for $n=1$. Then we compute:
\begin{eqnarray}
b^{\prec n} \prec b & =& (b\prec b^{\prec n-1})\prec b,\\
& =& b\prec (b^{\prec n-1}\prec b + b\prec b^{\prec n-1}),\\
& =& b\prec ((n-1)\ b^{\prec n}+  b^{\prec n}),\\
& =& n\ b\prec  b^{\prec n} = n\  b^{\prec n+1}).
\end{eqnarray}

Similarly, we compute by induction:
\begin{eqnarray}
(\cdots((a\prec \underbrace{b)\prec b) \cdots )\prec b}_{n} &=& (n-1)!\  (a\prec(b^{\prec n-1}))\prec b ,\\
&=& (n-1)! \ (a\prec(b^{\prec n-1}\prec b + b \prec b^{\prec n-1})),\\
&=& n! \ a\prec b^{\prec n}.
\end{eqnarray}
\hfill $\square$

It is known (cf.~\cite{L2}) that the tensor product of a Leibniz 
algebra with a Zinbiel algebra is a Lie algebra. Here we prove a finer 
result which is valid without characteristic hypothesis.

\begin{prop}\label{tensor} Over any field $k$ the tensor product $\g \otimes \mathfrak{R}$ of the 
Leibniz algebra $\g$ with the Zinbiel algebra $\mathfrak{R}$ is a pre-Lie algebra.
\end{prop}
\Proof For $x,y\in \g$ and $a,b\in \mathfrak{R}$ we define
$$\{x\otimes a, y\otimes b\} := [x,y]\otimes a\prec b .$$

We compute
$$\begin{array}{l}
\{\{x\otimes a,y\otimes b\},z\otimes c\}- \{x\otimes a,\{y\otimes 
b,z\otimes c\}\} \\
\quad =  [[x,y],z]\otimes (a\prec b)\prec c - [x,[y,z]] \otimes a\prec 
(b\prec c)\\
  \quad =  [[x,y],z]\otimes \big(a\prec (b\prec c)+  a\prec (c\prec 
b)\big) +\big(-[[x,y],z]+[[x,z],y]\big) \otimes a\prec (b\prec c) \\
  \quad =  [[x,y],z]\otimes  a\prec (c\prec b) +[[x,z],y] \otimes a\prec 
(b\prec c) .
\end{array}$$
As a consequence the associator of the binary operation $\{-,-\}$ is 
right-symmetric. We have proved that the operation $\{-,-\}$ is pre-Lie. 
\hfill $\square$

\medskip
 The category of pre-Lie algebras in prime characteristic have been studied by A. Dzhumadil'daev in \cite{DZ2}. In particular A. Dzumadil'daev, by theorem 1.1 in \cite{DZ2}, proves that there is a Jacobson formula for the $p$-th power of a sum of two elements of a pre-Lie algebra. Moreover, he introduces the notion of restricted pre-Lie algebra:

\subsection{Definition \cite{DZ2}}
A \emph{restricted pre-Lie} algebra is a pre-Lie algebra $A$ over $k$ such that $R_{a}^{p}=R_{a^{\{p\}}}$ for all $a\in A$.

The next theorem give us the Lie property of the Koszul duality in prime characteristic for the case of Leibniz and Zinbiel algebras.

\begin{thm}\label{pre-Lie}
Let $\g$ be a Leibniz algebra and $\mathfrak{R}$ a Zinbiel algebra. Then the tensor product $\g\otimes \mathfrak{R}$ is a restricted pre-Lie algebra with pre-Lie bracket given by:
\[ \,\{x\otimes a, y\otimes b\,\}:= \,[x,y]\otimes (a\prec b) \; \]
and $p$-map given by 
\[(y \otimes b)^{\{p\}}= [\cdots [[\underbrace{y,y],y]\cdots ,y}_{p}]\otimes (\cdots((\underbrace{b\prec b)\prec b) \cdots \prec b)}_{p}.\]
\end{thm}

\proof From Proposition \ref{tensor} it follows that the bracket $\{-\otimes -, -\otimes -\} $ defined on 
$\g\otimes \mathfrak{R}$ is indeed a pre-Lie bracket.

Moreover, we claim that  $\g\otimes \mathfrak{R}$ is a restricted pre-Lie algebra. Indeed we will prove that $R_{y\otimes b}^{p}=R_{(y\otimes b)^{\{p\}}}=0$. Since $\g$ is a Leibniz algebra, we have the relations $$[x, [y,z]]=-[x,[z,y]],\;\;\; [x,[y,y]]=0 $$ for all $x,y \in \g$. Using the relations above we get $\Big[x, [\cdots [[\underbrace{y,y],y]\cdots , y]}_{p}\Big]=0$. Thus,
 \[\, \{x \otimes a, (y\otimes b)^{ \{p\} }\}= \Big[x, [\cdots [[\underbrace{y,y],y]\cdots , y]}_{p}\Big] \otimes a \prec  \Big( (\cdots ((\underbrace{b\prec b)\prec b)\cdots)\prec b}_{p}\Big)=0. \]
Besides, by Lemma \ref{zinb} above we obtain that $R_{(y\otimes b)}^{p}=0.$ Therefore $\g\otimes \mathfrak{R}$ is a restricted pre-Lie algebra as claimed.
\hfill $\square$

\begin{cor}
Let $\g$ be a Leibniz algebra and $\mathfrak{R}$ a Zinbiel algebra. Then the tensor product $\g\otimes \mathfrak{R}$ is a restricted Lie algebra with Lie bracket given by:
\[ \,[x\otimes a, y\otimes b\,]:= \,[x,y\,]\otimes (a\prec b)- [\,y,x\,]\otimes (b \prec a)\; \]
and $p$-map given by \[(y \otimes b)^{[p]}= [\cdots [[\underbrace{y,y],y]\cdots ,y}_{p}]\otimes (\cdots((\underbrace{b\prec b)\prec b) \cdots \prec b)}_{p}.\]
\end{cor}
\proof
By Corollary 2.4 in \cite{DZ2} and Theorem \ref{pre-Lie} the corresponding Lie algebra $(\g\otimes \mathfrak{R})_{Lie}$ is a restricted Lie algebra with $p$-map given by $(y\otimes b) \mapsto (y\otimes b)^{ \{p\} }$ and the corollary is proved.
\hfill $\square$

\end{document}